\documentclass[12pt]{article}
\usepackage{latexsym,amsmath}
\usepackage{amssymb}
\usepackage[margin=1.25in]{geometry}

\newtheorem{thm}{Theorem}[section]

\newtheorem{lemma}[thm]{Lemma}

\newtheorem{conj}[thm]{Conjecture}

\newcommand{\Var}{{\operatorname{\rm Var}}}

\newcommand{\beq}[1]{\begin{equation}\label{#1}}
\newcommand{\enq}[0]{\end{equation}}

\newcommand{\qed}[0]{{\hspace*{\fill}\mbox{$\Box$}}}

\newcommand{\ul}[0]{\underline}

\newcommand{\cA}[0]{{\cal A}}
\newcommand{\cB}[0]{{\cal B}}

\newcommand{\cE}[0]{{\cal E}}

\newcommand{\cH}[0]{{\cal H}}
\newcommand{\cI}[0]{{\cal I}}

\newcommand{\cL}[0]{{\cal L}}

\newcommand{\cO}[0]{{\cal O}}

\newcommand{\N}{{\mathbb N}}

\newcommand{\R}{{\mathbb R}}

\newcommand{\gb}[0]{\beta}

\newcommand{\gd}[0]{\delta}

\newcommand{\gl}[0]{\lambda}
\newcommand{\gL}[0]{\Lambda}

\newcommand{\gs}[0]{\sigma}

\title{Bounding the partition function of spin-systems}

\author{David J. Galvin\thanks{Department of Mathematics, University of Pennsylvania, 209 South 33th Street, Philadelphia PA 19104. This work was begun while the author
was a member of the Institute for Advanced Study, Einstein Drive,
Princeton NJ 08540 and was supported in part by NSF grant
DMS-0111298.}}
\date{Appeared 2006}
\begin{document}
\maketitle

\begin{abstract}
With a graph $G=(V,E)$ we associate a collection of non-negative
real weights $\cup_{v\in V}\{\gl_{i,v}:1\leq i \leq m\} \cup
\cup_{uv \in E} \{\gl_{ij,uv}:1\leq i \leq j \leq m\}$. We consider
the probability distribution on $\{f:V\rightarrow\{1,\ldots,m\}\}$
in which each $f$ occurs with probability proportional to $\prod_{v
\in V}\gl_{f(v),v}\prod_{uv \in E}\gl_{f(u)f(v),uv}$. Many
well-known statistical physics models, including the Ising model
with an external field and the hard-core model with non-uniform
activities, can be framed as such a distribution. We obtain an upper
bound, independent of $G$, for the partition function (the
normalizing constant which turns the assignment of weights on
$\{f:V\rightarrow\{1,\ldots,m\}\}$ into a probability distribution)
in the case when $G$ is a regular bipartite graph. This generalizes
a bound obtained by Galvin and Tetali who considered the simpler
weight collection $\{\gl_i:1 \leq i \leq m\} \cup \{\gl_{ij}:1 \leq
i \leq j \leq m\}$ with each $\gl_{ij}$ either $0$ or $1$ and with
each $f$ chosen with probability proportional to $\prod_{v \in
V}\gl_{f(v)}\prod_{uv \in E}\gl_{f(u)f(v)}$. Our main tools are a
generalization to list homomorphisms of a result of Galvin and
Tetali on graph homomorphisms and a straightforward second-moment
computation.

\end{abstract}

\newpage

\section{Introduction and statement of results}

Let $G=(V(G),E(G))$ and $H=(V(H),E(H))$ be non-empty graphs. Set
$$
Hom(G,H)=\{f:V(G)\rightarrow V(H): uv \in E(G)\Rightarrow f(u)f(v)
\in E(H)\}
$$
(that is, $Hom(G,H)$ is the set of graph homomorphisms from $G$ to
$H$). In \cite{GalvinTetali}, the following result is obtained
using entropy considerations (and in particular, Shearer's Lemma).
\begin{thm} \label{thm-from.gt.unweighted}
For any graph $H$ and any $d$-regular $N$-vertex bipartite graph
$G$,
$$
|Hom(G,H)| \leq |Hom(K_{d,d},H)|^\frac{N}{2d}
$$
where $K_{d,d}$ is the complete bipartite graph with $d$ vertices
in each partition class.
\end{thm}

In \cite{GalvinTetali} Theorem \ref{thm-from.gt.unweighted} is
extended to a result on weighted graph homomorphisms. To each
$i\in V(H)$ assign a positive pair of weights $(\gl_i,\mu_i)$.
Write $(\gL,\rm{M})$ for the set of weights. For a bipartite graph
$G$ with bipartition classes $\cE_G$ and $\cO_G$ give each $f \in
Hom(G,H)$ weight
$$
w^{(\gL,\rm{M})}(f) :=\prod_{v\in \cE_G}\gl_{f(v)} \prod_{v\in
\cO_G}\mu_{f(v)}.
$$
The constant that turns this
assignment of weights on $Hom(G,H)$ into a probability
distribution is
$$
Z^{(\gL,\rm{M})}(G,H):=\sum_{f\in Hom(G,H)}w^{(\gL,\rm{M})}(f).
$$
The following is proved in \cite{GalvinTetali}.
\begin{thm} \label{thm-from.gt.most.general}
For any graph $H$, any set $(\gL,\rm{M})$ of positive weights on
$V(H)$ and any $d$-regular $N$-vertex bipartite graph $G$,
$$
Z^{(\gL,\rm{M})}(G,H) \leq
\left(Z^{(\gL,\rm{M})}(K_{d,d},H)\right)^{\frac{N}{2d}}.
$$
\end{thm}
Taking all weights to be $1$, Theorem
\ref{thm-from.gt.most.general} reduces to Theorem
\ref{thm-from.gt.unweighted}.

\medskip

In this note we consider a more general weighted model. Fix $m \in
\N$ and a graph $G=(V, E)$. With each $1 \leq i \leq m$ and $v \in
V$ associate a non-negative real weight $\gl_{i,v}$ and with each $1
\leq i \leq j \leq m$ and $uv \in E$ associate a non-negative real
weight $\gl_{ij,uv}$. Set $\gl_{ij,uv}:=\gl_{ji,uv}$ for $i > j$.
Write $W$ for the collection of weights and for each $f:V
\rightarrow \{1,\ldots,m\}$ set
$$
w^W(f) = \prod_{v \in V} \gl_{f(v),v}\prod_{uv\in
E}\gl_{f(u)f(v),uv}
$$
and
$$
Z^W(G)=\sum_{f:V \rightarrow \{1,\ldots,m\}} w^W(f).
$$

We may put all this in the framework of a well-known mathematical
model of physical spin systems. We think of the vertices of $G$ as
particles and the edges as bonds between pairs of particles
(typically a bond represents spatial proximity), and we think of
$\{1,\ldots, m\}$ as the set of possible spins that a particle may
take. For each $v \in V$ we think of the weights $\gl_{\cdot,v}$ as
a measure of the likelihood of seeing the different spins at $v$;
furthermore, for each $uv \in E$ we think of the weights
$\gl_{\cdot,uv}$ as a measure of the likelihood of seeing the
different spin-pairs across the edge $uv$. The probability of a
particular spin configuration is thus proportional to the product
over the vertices of $G$ of the weights of the spins times the
product over the edges of $G$ of the weights of the spin-pairs. In
this language $Z^W(G)$ is the partition function of the model
--- the normalizing constant that turns the above-described system
of weights on the set of spin configurations into a probability
measure.

An example of such a model is the {\em hard-core} (or independent
set) model. Here $m=2$ and the system of weights $W_{hc}$ is given
by
$$
\begin{array}{ccc}
\gl_{i,v} = \left\{
\begin{array}{ll}
\gl_v & \mbox{if $i=1$}\\
1 & \mbox{if $i=2$}
\end{array}
\right. & ~~~\mbox{and}~~~ & \gl_{ij,uv} = \left\{
\begin{array}{ll}
0 & \mbox{if $i=j=1$}\\
1 & \mbox{otherwise,}
\end{array}
\right.
\end{array}
$$
and so
$$
Z^{W_{hc}}(G)=\sum_{f:V\rightarrow \{1,2\}} \left(\prod_{v:f(v)=1}
\gl_v\right)\left({\bf 1}_{\{\not \exists uv \in
E:f(u)=f(v)=1\}}\right) = \sum_{I \in \cI(G)} \prod_{v \in I} \gl_v
$$
is a weighted sum of independent sets in $G$. (Recall that $I
\subseteq V$ is {\em independent} in $G$ if for all $u,v \in I$, $uv
\not \in E$. We write $\cI(G)$ for the collection of independent
sets in $G$.)

The hard-core model is a {\em hard-constraint} model in which all of
the edge-weights are either $0$ or $1$, and the r\^{o}le of these
weights is to exclude certain configurations from contributing to
the partition function. We now consider the best known example of a
{\em soft-constraint} model (one in which all configurations are
potentially allowable), the {\em Ising} model. Here $m=2$ and there
are two parameters, $\gb,h \in \R$. It is convenient to take the set
of spins to be $\{+1,-1\}$. The system $W_{Ising,\gb,h}$ of weights
on $\{+1, -1\}$ is as follows.
\begin{eqnarray*}
\gl_{+1,v} & = & e^h ~\mbox{for all $v \in V$} \\
\gl_{-1,v} & = & e^{-h} ~\mbox{for all $v \in V$} \\
\gl_{ii,uv} & = & e^{-\gb} ~\mbox{for $i\in \{+1,-1\}$ and all $uv \in E$ and} \\
\gl_{+1-1,uv} & = & e^{\gb} ~\mbox{for all $uv \in E$.}
\end{eqnarray*}
For each $\gs:V\rightarrow \{+1,-1\}$ we have
$$
w^{W_{Ising,\gb,h}}(\gs)=\exp\left\{-\gb\sum_{uv \in E}
\gs(u)\gs(v)+h\sum_{v \in V(G} \gs(v)\right\}.
$$
Then $Z^{W_{Ising,\gb,h}}(G)=\sum_\gs w^{W_{Ising,\gb,h}}(\gs)$ is
the partition function of the Ising model on $G$ with inverse
temperature $|\gb|$ and external field $h$. (If $\gb > 0$, we are in
the anti-ferromagnetic case, where configurations with many
$+1$-$-1$ edges are favoured; if $\gb < 0$, we are in the
ferromagnetic case, where configurations with few $+1$-$-1$ edges
are favoured.)

\medskip

Let us now set up the notation for our main result. For
completeness, we choose to make the straightforward generalization
from regular bipartite graphs to $(a,b)$-biregular graphs, that is,
bipartite graphs in which one partition class, which we shall label
$\cE_G$, consists of vertices of degree $a$ and the other class,
$\cO_G$, consists of vertices of degree $b$. For $v \in \cO_G$ write
$\{n_1(v),\ldots,n_b(v)\}$ for the set of neighbours of $v$.

Let $G$ be such a graph and let $W$ be a collection of weights on
$G$. Give labels $w_1,\ldots,w_b$ to the degree $a$ vertices of
$K_{a,b}$ (the complete bipartite graph with $a$ vertices in one
partition class and $b$ in the other) and labels $z_1,\ldots,z_a$ to
the degree $b$ vertices. For $v \in \cO_G$ write $W^v$ for the
following system of weights on $K_{a,b}$:
$$
\gl^v_{i,u} = \left\{
\begin{array}{ll}
\gl_{i,v} & \mbox{if $u=z_k$ for some $1\leq k \leq a$} \\
\gl_{i,n_k(v)} & \mbox{if $u =w_k$ for some $1\leq k \leq b$}
\end{array}
\right.
$$
and, for $1\leq k \leq b$ and $1\leq \ell \leq a$, $\gl_{ij,w_k
z_\ell} = \gl_{ij,n_k(v)v}$. Our main result is a generalization of
Theorem \ref{thm-from.gt.most.general} to the following.
\begin{thm} \label{thm-most.general.weighted}
For any $(a,b)$-biregular graph $G$ and any system of weights $W$,
$$
Z^{W}(G) \leq \prod_{v \in {\cal O}_G}
\left(Z^{W^v}(K_{a,b})\right)^{\frac{1}{a}}.
$$
\end{thm}
Taking $a=b=d$,
$$
\gl_{i,v} =
\begin{array}{ccc}
\left\{
\begin{array}{ll}
\gl_i & \mbox{if $v \in \cE_G$} \\
\mu_i & \mbox{if $v \in \cO_G$}
\end{array}
\right. & ~~~\mbox{and}~~~ & \gl_{ij,uv} = \left\{
\begin{array}{ll}
1 & \mbox{if $ij \in E(H)$} \\
0 & \mbox{otherwise,}
\end{array}
\right.
\end{array}
$$
Theorem \ref{thm-most.general.weighted} reduces to Theorem
\ref{thm-from.gt.most.general}.

\medskip

Let us consider an application of Theorem
\ref{thm-most.general.weighted} to the antiferromagnetic ($\gb > 0$)
Ising model without external field ($h=0$) on a $d$-regular,
$N$-vertex bipartite graph $G$. A trivial lower bound on
$Z^{W_{Ising,\gb,h}}(G)$,
\begin{equation} \label{inq-lb.on.ising}
\exp\left\{\frac{\gb dN}{2}\right\} \leq Z^{W_{Ising,\gb,h}}(G),
\end{equation}
is obtained by considering the configuration in which one partition
class of $G$ is mapped entirely to $+1$ and the other entirely to
$-1$. Applying Theorem \ref{thm-most.general.weighted} we obtain as
an upper bound
\begin{eqnarray}
Z^{W_{Ising,\gb,h}}(G) & \leq &
Z^{W_{Ising,\gb,h}}(K_{d,d})^\frac{N}{2d} \nonumber \\
& \leq & \left(2^{2d}e^{\gb d^2}\right)^\frac{N}{2d} \label{ub.exp} \\
& = & 2^N\exp\left\{\frac{\gb dN}{2}\right\}
\label{inq-ub.on.ising}.
\end{eqnarray}
In (\ref{ub.exp}) we are using that there are $2^{2d}$ possible
configurations on $K_{d,d}$ and that each has weight at most $e^{\gb
d^2}$. Combining (\ref{inq-lb.on.ising}) and (\ref{inq-ub.on.ising})
we obtain the following bounds on the {\em free-energy} of the Ising
model, the quantity
$F^{W_{Ising,\gb,h}}(G):=\log(Z^{W_{Ising,\gb,h}}(G))/N$:
$$
\frac{\gb d}{2} \leq F^{W_{Ising,\gb,h}}(G) \leq \frac{\gb d}{2}
+\ln 2.
$$
Note that these bounds are absolute (independent of $G$ and $N$),
and asymptotically tight in the case of a family of graphs
satisfying $\gb d = \omega(1)$.

\medskip

We give the proof of Theorem \ref{thm-most.general.weighted} in
Section \ref{sec-main.proof}. An important tool in the proof is an
extension of Theorem \ref{thm-from.gt.unweighted} to the case of
list homomorphisms, which we now discuss. Let $H$ and $G$ be
non-empty graphs. To each $v \in V(G)$ associate a set $L(v)
\subseteq V(H)$ and write $\cL(G,H)$ for $\{L(v):v \in V(G)\}$. A
{\em list homomorphism from $G$ to $H$ with list set $\cL(G,H)$} is
a homomorphism $f \in Hom(G,H)$ satisfying $f(v) \in L(v)$ for all
$v \in V(G)$. Write $Hom^{\cL(G,H)}(G,H)$ for the set of all list
homomorphism from $G$ to $H$ with list set $\cL(G,H)$.

The notion of a list homomorphism is a generalization of that of a
homomorphism. Indeed, if $L(v)=V(H)$ for all $v \in V(G)$ then
$Hom^{\cL(G,H)}(G,H)$ is the same as $Hom(G,H)$. List homomorphisms
also generalize the well-studied notion of list colourings of a
graph (see {\em e.g.} \cite[Chapter 5]{Diestel} for an
introduction). Recall that if a list $L(v)$ of potential colours is
assigned to each vertex $v$ of a graph $G$, then a list colouring of
$G$ (with list set $\cL(G)=\{L(v):v \in V(G)\}$) is a function
$\chi:V(G) \rightarrow \cup_{v \in V(G)} L(v)$ satisfying the
property that $\chi$ is a proper colouring ({\em i.e.}, $\chi(u)\neq
\chi(v)$ for all $uv \in E(G)$) that respects the lists ({\em i.e.},
$\chi(v) \in L(v)$ for all $v \in V(G)$). The set of list colourings
of $G$ with list set $\cL(G)$ is exactly the set
$Hom^{\cL(G)}(G,H_{\cL(G)})$ where $H_{\cL(G)}$ is the complete
loopless graph on vertex set $\cup_{v \in V(G)} L(v)$.

In the discussion that follows we fix an $(a,b)$-biregular graph
$G$. We also fix $H$ and $\cL(G,H)$ and for convenience of notation
we often suppress dependence on $G$ and $H$. For $v \in \cO_G$ write
$\cL^v$ for the list set on $K_{a,b}$ in which each vertex of degree
$b$ gets list $L(v)$ and the vertices of degree $a$ get the lists
$L(n_1(v)),\ldots,L(n_b(v))$ (each one occurring exactly once) where
$\{n_1(v),\ldots,n_b(v)\}$ is the set of neighbours of $v$. (Recall
that $K_{a,b}$ is the complete bipartite graph with $a$ vertices in
one partition class and $b$ in the other.) We generalize Theorem
\ref{thm-from.gt.unweighted} to the following result, whose proof is
given in Section \ref{sec-list.proof}.
\begin{thm} \label{thm-unweighted.list}
For any graph $H$, any $(a,b)$-biregular graph $G$ and any list set
$\cL$,
$$
\left|Hom^\cL(G,H)\right|\leq
    \prod_{v \in {\cal O}_G}
    \left(\left|Hom^{\cL^v}(K_{a,b},H)\right|\right)^{\frac{1}{a}}.
$$
\end{thm}
Taking $a=b=d$ and $L(v)=V(H)$ for all $v \in V(G)$, Theorem
\ref{thm-unweighted.list} reduces to Theorem
\ref{thm-from.gt.unweighted}.

\medskip

Before turning to proofs, we pause to make a conjecture. The point
of departure for this note and for \cite{GalvinTetali} is a result
of Kahn \cite{Kahn} bounding the number of independent sets in a
$d$-regular, $N$-vertex bipartite graph $G$ by
\begin{equation} \label{ind}
|\cI(G)| \leq |\cI(K_{d,d})|^\frac{N}{2d}.
\end{equation}
Kahn conjectured in \cite{Kahn} that for an arbitrary graph $G$ it
should hold that
\begin{equation} \label{ind.conj}
|\cI(G)| \leq \prod_{uv \in
E(G)}|\cI(K_{d(u),d(v)})|^\frac{1}{d(u)d(v)}.
\end{equation}
where $d(u)$ denotes the number of neighbours of $u$ in $G$. Note
that (\ref{ind}) is a special case of (\ref{ind.conj}), and that
(\ref{ind.conj}), if true, would be tight for any $G$ which is the
union of complete bipartite graphs.

At the moment we see no reason not to conjecture the following,
which stands in relation to Theorem \ref{thm-most.general.weighted}
as (\ref{ind.conj}) does to (\ref{ind}).
\begin{conj} \label{conj-most.general}
Let $G$ be any graph and $W$ any collection of weights on $G$. For
each $u \in V(G)$ let $\{n_1(u),\ldots,n_{d(u)}(u)\}$ be the set of
neighbours of $u$. For each edge $uv \in E(G)$, label the degree
$d(u)$ vertices of $K_{d(u),d(v)}$ by
$w_1(u,v),\ldots,w_{d(v)}(u,v)$ and the degree $d(v)$ vertices by
$z_1(u,v),\ldots,z_{d(u)}(u,v)$. Let $W^{uv}$ be the collection of
weights on $K_{d(u),d(v)}$ given by
$$
\gl^{u,v}_{i,w_j(u,v)} = \gl_{i,n_j(v)},~~ \gl^{u,v}_{i,z_j(u,v)} =
\gl_{i,n_j(u)} ~~~\mbox{and}~~~ \gl^{u,v}_{ij,w_j(u,v)z_k(u,v)} =
\gl_{ij,n_j(v)n_k(u)}.
$$
Then
$$
Z^W(G) \leq \prod_{uv \in E(G)}
Z^{W^{uv}}(K_{d(u),d(v)})^{\frac{1}{d(u)d(v)}}.
$$
\end{conj}
Exactly as Theorem \ref{thm-most.general.weighted} follows from
Theorem \ref{thm-unweighted.list} (as will be described in Section
\ref{sec-main.proof}), Conjecture \ref{conj-most.general} would
follow from the following conjecture concerning list homomorphisms.
\begin{conj} \label{conj-most.general.list}
Let $G$ and $H$ be any graphs and $\cL$ any list set. Let $\cL^{uv}$
be the list set on $K_{d(u),d(v)}$ given by
$$
L^{u,v}(w_j(u,v)) = L(n_j(v)) ~~~\mbox{and}~~~ L^{u,v}(z_j(u,v)) =
L(n_j(u))
$$
(with the notation as in Conjecture \ref{conj-most.general}). Then
$$
|Hom^\cL(G,H)| \leq \prod_{uv \in E(G)}
|Hom^{\cL^{uv}}(K_{d(u),d(v)})|^{\frac{1}{d(u)d(v)}}.
$$
\end{conj}

\section{Proof of Theorem \ref{thm-unweighted.list}}
\label{sec-list.proof}

We derive Theorem \ref{thm-unweighted.list} from the following more
general statement.
\begin{thm} \label{thm-most.general}
Let $G$ be a bipartite graph with partition classes $\cE_G$ and
$\cO_G$, $H$ an arbitrary graph and $\cL=\cL(G,H)$ a list set.
Suppose that there is $m$, $t_1$ and $t_2$ and families $\cA =
\{A_i:1 \leq i \leq m\}$ and $\cB = \{B_i:1 \leq i \leq m\}$ with
each $A_i \subseteq \cE_G$ and each $B_i \subseteq \cO_G$ such that
each $v \in \cE_G$ is contained in at least $t_1$ members of $\cA$
and each $u \in \cO_G$ is contained in at least $t_2$ members of
$\cB$. Then
$$
|Hom^\cL(G,H)| \leq \prod_{i=1}^m \left(\sum_{\ul{x}\in \prod_{v
\in A_i} L(v)} |C^{\ul{x}}(A_i,B_i)|^{\frac{t_1}{t_2}}
\right)^{\frac{1}{t_1}}
$$
where, for each $1 \leq i \leq m$ and each $\ul{x} \in \prod_{v
\in A_i} L(v)$,
$$
C^{\ul{x}}(A_i,B_i)=\left\{f:B_i \rightarrow V(H):\begin{array}{l}
\forall ~v \in B_i,~f(v)\in L(v)~\mbox{and} \\
\forall ~v \in B_i,~u \in A_i ~\mbox{with $uv \in
E(G)$},~(\ul{x})_uf(v) \in E(H)
\end{array}\right\}
$$
is the set of extensions of the partial list homomorphism $\ul{x}$
on $A_i$ to a partial list homomorphism on $A_i \cup B_i$.
\end{thm}
To obtain Theorem \ref{thm-unweighted.list} from Theorem
\ref{thm-most.general} we take $\cA=\{N(v):v \in \cO_G\}$ and
$\cB=\{\{v\}:v \in \cO_G\}$ where $N(v)=\{u \in V(G):uv \in
E(G)\}$ so that $t_1=a$ and $t_2=1$, and note that in this case
$\sum_{\ul{x}\in \prod_{u \in A_v} L(u)}
|C^{\ul{x}}(N(v),\{v\})|^a$ is precisely
$|Hom^{\cL^v}(K_{a,b},H)|$.

\medskip

The proof of Theorem \ref{thm-most.general} uses entropy
considerations, which for completeness we very briefly review here.
Our treatment is mostly copied from \cite{Kahn}. For a more thorough
discussion, see {\em e.g.} \cite{McE}. In what follows ${\bf X}$,
${\bf Y}$ etc. are discrete random variables, which in our usage are
allowed to take values in any finite set.

The {\em entropy} of ${\bf X}$ is
$$\
H({\bf X}) = \sum_x p(x)\log\frac{1}{p(x)},
$$
where we write $p(x)$ for ${\mathbb P}({\bf X}=x)$ (and extend this
convention in natural ways below). The {\em conditional entropy} of
${\bf X}$ given ${\bf Y}$ is
$$
H({\bf X}|{\bf Y}) ={\mathbb E} H({\bf X}|\{{\bf Y}=y\})
=\sum_yp(y)\sum_xp(x|y)\log\frac{1}{p(x|y)}.
$$
Notice that we are also writing $H({\bf X}|Q)$ with $Q$ an event (in
this case $Q=\{{\bf Y} =y\}$):
$$
H({\bf X}|Q)=\sum p(x|Q)\log\frac{1}{p(x|Q)}.
$$
We have the inequalities
$$
H({\bf X}) \leq \log |{\rm range}({\bf
X})|~~~~ \mbox{(with equality if ${\bf X}$ is uniform),}
$$
$$
H({\bf X}|{\bf Y}) \leq H({\bf X}),
$$
and more generally,
\begin{equation} \label{dropping}
\mbox{if ${\bf Y}$ determines ${\bf Z}$ then $H({\bf X}|{\bf Y})
\leq H({\bf X}|{\bf Z})$.}
\end{equation}
For a random vector ${\bf X}=({\bf X}_1,\ldots ,{\bf X}_n)$ there is
a chain rule
\begin{equation} \label{chain}
H({\bf X}) =H({\bf X}_1)+H({\bf X}_2|{\bf X}_1)+\cdots
    +H({\bf X}_n|{\bf X}_1,\ldots, {\bf X}_{n-1}).
\end{equation}
Note that (\ref{dropping}) and (\ref{chain}) imply
\begin{equation} \label{sub}
H({\bf X}_1,\ldots, {\bf X}_n)\leq \sum H({\bf X}_i)
\end{equation}
We also have a conditional version of (\ref{sub}):
$$
H({\bf X}_1,\ldots, {\bf X}_n|{\bf Y})\leq \sum H({\bf X}_i|{\bf Y}).
$$
Finally we use a lemma of Shearer (see \cite[p. 33]{CFGS}). For a
random vector ${\bf X}=({\bf X}_1,\ldots, {\bf X}_m)$ and
$A\subseteq \{1,\ldots,m\}$, set ${\bf X}_A=({\bf X}_i:i\in A)$.
\begin{lemma}
\label{Lshearer} Let ${\bf X}=({\bf X}_1,\ldots, {\bf X}_m)$ be a
random vector, ${\bf Y}$ a random variable and ${\cal A}$ a
collection of subsets (possibly with repeats) of $\{1,\ldots,m\}$,
with each element of $\{1,\ldots,m\}$ contained in at least $t$
members of ${\cal A}$. Then
$$
H({\bf X})\leq
    \frac{1}{t}\sum_{A\in{\cal A}}H({\bf X}_A) ~~~~~~\mbox{and}~~~~~~H({\bf X}|{\bf Y})\leq
    \frac{1}{t}\sum_{A\in{\cal A}}H({\bf X}_A|{\bf Y}).
$$
\end{lemma}

\noindent {\em Proof of Theorem \ref{thm-most.general}: }We follow
closely the proof of \cite[Lemma 3.1]{GalvinTetali}. Let ${\bf f}$
be a uniformly chosen member of $Hom^\cL(G,H)$. For each $1 \leq i
\leq m$ and each $\ul{x} \in \prod_{v \in A_i} L(v)$ let
$p_i(\ul{x})$ be the probability that ${\bf f}$ restricted to $A_i$
is $\ul{x}$. With the key inequalities justified below (the
remaining steps follow in a straightforward way from the properties
of entropy just established) we have
\begin{eqnarray}
H({\bf f}) & = & H({\bf f}|_{\cE_G})+H({\bf f}|_{\cO_G}~|~{\bf
f}|_{\cE_G})
\nonumber \\
& \leq & \frac{1}{t_1}\sum_{i=1}^m H({\bf f}|_{A_i}) +
\frac{1}{t_2}\sum_{i=1}^m H({\bf f}|_{B_i}~|~{\bf f}|_{\cE_G})
\label{shearer.3}\\
& \leq & \frac{1}{t_1}\sum_{i=1}^m \left(H({\bf
f}|_{A_i})+\frac{t_1}{t_2}H({\bf f}|_{B_i}~|~{\bf
f}|_{A_i})\right) \label{inq_vs_eq} \\
& = & \frac{1}{t_1}\sum_{i=1}^m \sum_{\ul{x}\in \prod_{v \in A_i}
L(v)}
         \left(p_i(\ul{x})\log \frac{1}{p_i(\ul{x})}
         + \frac{t_1}{t_2}p_i(\ul{x})H({\bf f}|_{B_i}~|~\{{\bf f}|_{A_i}=\ul{x}\})\right)
         \nonumber \\
& \leq & \frac{1}{t_1}\sum_{i=1}^m \sum_{\ul{x}\in \prod_{v \in A_i}
L(v)}
p_i(\ul{x})\log\frac{|C^{\ul{x}}(A_i,B_i)|^{\frac{t_1}{t_2}}}{p_i(\ul{x})}
\nonumber \\
& \leq & \frac{1}{t_1}\sum_{i=1}^m \log
\left(\sum_{\ul{x}\in \prod_{v \in A_i} L(v)}
|C^{\ul{x}}(A_i,B_i)|^{\frac{t_1}{t_2}}\right). \label{jensen.3}
\end{eqnarray}
In (\ref{shearer.3}) we use Shearer's Lemma twice, once with $\cA$
as the covering family and once with $\cB$, and in (\ref{jensen.3})
we use Jensen's inequality. In (\ref{inq_vs_eq}) we would have
equality if it happened that for each $i$, $A_i$ included all the
neighbours of $B_i$, since ${\bf f}|_{B_i}$ depends only on the
values of ${\bf f}$ on $B_i$'s neighbours. It is easy, however, to
construct examples where $H({\bf f}|_{B_i}~|~{\bf f}|_{\cE_G}) <
H({\bf f}|_{B_i}~|~{\bf f}|_{A_i})$ when $A_i$ does not include all
the neighbours of $B_i$.

The theorem now follows from the equality $H({\bf f}) =
\log|Hom^\cL(G,H)|$. \qed

\section{Proof of Theorem \ref{thm-most.general.weighted}}
\label{sec-main.proof}

By continuity we may assume that all weights are rational and
non-zero. By scaling appropriately we may also assume that $0 <
\gl_{ij,uv} \leq 1$ for all $i$, $j$ and $uv \in E(G)$ (we will
later think of the $\gl_{ij,uv}$'s as probabilities).

Set $N=|V(G)|$ and
$$
\gl_{vmin}= \min_{i,w} \gl_{i,w}, ~~~\gl_{vmax}= \max_{i,w}
\gl_{i,w}~~~\mbox{and}~~~\gl_{emin}=\min_{ij,vw} \gl_{ij,vw}.
$$
Also, set $w_{min}=\gl_{vmin}^N \gl_{emin}^{abN/(a+b)}$; this is a
lower bound on $w^W(f)$ for all $f:V(G)\rightarrow \{1,\ldots, m\}$
(observe that an $(a,b)$-biregular graph $G$ on $N$ vertices has
$|\cE_G|=bN/(a+b)$, $|\cO_G|=aN/(a+b)$ and $|E(G)|=abN/(a+b)$) as
well as a lower bound on $w^{W^v}(f)$ for all $v \in \cO_G$ and all
$f:V(K_{a,b})\rightarrow \{1,\ldots, m\}$.

Choose $C\geq 1$ large enough that $C\gl_{i,v} \in \N$ for all $1
\leq i \leq m$ and $v \in V(G)$. For each $i$ and $v$ let
$S_{i,v}$ be a set of size $C\gl_{i,v}$, with all the $S_{i,v}$'s
disjoint. Let $H$ be the graph on vertex set $\cup_{i,v} S_{i,v}$
with $xy \in E(H)$ iff $x \in S_{i,v}$ and $y \in S_{j,w}$ for
some $i,j,v,w$ with $vw \in E(G)$. For each $v \in V(G)$ let
$L(v)=\cup_i S_{i,v}$ and set $\cL=\{L(v): v \in V(G)\}$. For each
$g:V(G) \rightarrow \{1,\ldots, m\}$ and each subgraph
$\widetilde{H}$ of $H$ (on the same vertex set as $H$) set
$$
\cH_g(G,\widetilde{H}) = \{f \in Hom^\cL(G,\widetilde{H}):f(v) \in
S_{g(v),v} \mbox{ for all $v \in V(G)$}\}.
$$
Note that $\cH_g(G,H)$ is exactly $\{f:V(G)\rightarrow V(H):f(v)
\in S_{g(v),v} \mbox{ for all $v \in V(G)$}\}$ and so
$|\cH_g(G,H)|=C^N \prod_{v \in V(G)} \gl_{g(v),v}$. Note also that
for $g \neq g'$ we have $\cH_g(G,\widetilde{H}) \cap
\cH_{g'}(G,\widetilde{H})=\emptyset$ and that $Hom^\cL(G,
\widetilde{H}) = \cup_g \cH_g(G,\widetilde{H})$.

For each $v \in \cO_G$, each $g:V(K_{a,b})\rightarrow
\{1,\ldots,m\}$ and each $\widetilde{H}$, set
$$
\cH^v_g(K_{a,b},\widetilde{H})=\left\{f \in
Hom^{\cL^v}(K_{a,b},\widetilde{H}):\begin{array}{ll} f(w_k) \in
S_{g(w_k),n_k(v)}, & 1 \leq k \leq b \\
f(z_k) \in S_{g(z_k),v}, & 1 \leq k \leq a \end{array}\right\},
$$
where the notation is as established before the statements of
Theorems \ref{thm-most.general.weighted} and
\ref{thm-unweighted.list}. Note that for $g \neq g'$ we have
$\cH^v_g(K_{a,b},\widetilde{H}) \cap
\cH^v_{g'}(K_{a,b},\widetilde{H})=\emptyset$ and that
$Hom^{\cL^v}(K_{a,b}, \widetilde{H}) = \cup_g
\cH^v_g(K_{a,b},\widetilde{H})$.

We will exhibit a subgraph $\widetilde{H}$ of $H$ which satisfies
\begin{equation} \label{inq-sec.mom.main}
\left|C^Nw^W(g)-|\cH_g(G,\widetilde{H})|\right| \leq
\gd(C)|\cH_g(G,\widetilde{H})|
\end{equation}
for all $g:V(G)\rightarrow \{1,\ldots,m\}$ and
\begin{equation} \label{inq-sec.mom.main.2}
\left|C^{a+b}w^{W^v}(g)-|\cH^v_g(K_{a,b},\widetilde{H})|\right| \leq
\gd(C)|\cH^v_g(K_{a,b},\widetilde{H})|
\end{equation}
for all $v \in \cO_G$ and $g:V(K_{a,b})\rightarrow
\{1,\ldots,m\}$, where $\gd(C)$ depends also on $N$, $a$, $b$ and
$W$ and tends to $0$ as $C$ tends to infinity (with $N$, $a$, $b$
and $W$ fixed). This suffices to prove the theorem, for we have
\begin{eqnarray}
\left|C^NZ^W(G)-|Hom^\cL(G, \widetilde{H})|\right| & \leq & \sum_g
\left|C^Nw^W(g)-|\cH_g(G,\widetilde{H})|\right|
\nonumber \\
& \leq & \gd(C)\sum_g
|\cH_g(G,\widetilde{H})| \nonumber \\
& = & \gd(C) |Hom^\cL(G, \widetilde{H})| \nonumber
\end{eqnarray}
and similarly, for each $v \in \cO_G$,
\begin{equation} \label{inq.2}
\left|C^{a+b}Z^{W^v}(K_{a,b})-|Hom^{\cL^v}(K_{a,b},
\widetilde{H})|\right| \leq \gd(C) |Hom^{\cL^v}(K_{a,b},
\widetilde{H})|
\end{equation}
and so
\begin{eqnarray}
C^NZ^W(G) & \leq & (1+\delta(C))|Hom^\cL(G,\widetilde{H})|
\nonumber \\
& \leq & (1+\delta(C))\prod_{v \in \cO_G}
\left(|Hom^{\cL^v}(K_{a,b},\widetilde{H})|\right)^{\frac{1}{a}}
\label{using.main.thm} \\
& \leq & C^N\frac{1+\delta(C)}{(1-\delta(C))^{\frac{N}{a+b}}}
\prod_{v \in \cO_G} \left(Z^{W^v}(K_{a,b})\right)^{\frac{1}{a}}
\label{using1}.
\end{eqnarray}
In (\ref{using.main.thm}) we use Theorem \ref{thm-unweighted.list}
while in (\ref{using1}) we use (\ref{inq.2}). Theorem
\ref{thm-most.general.weighted} follows since the constant in front
of the product in (\ref{using1}) can be made arbitrarily close to
$1$ (with $N$, $a$, $b$ and $W$ fixed) by choosing $C$ sufficiently
large.

The graph $\widetilde{H}$ will be a random graph defined as follows.
For each $xy \in E(H)$ with $x \in S_{i,v}$ and $y \in S_{j,w}$ we
put $xy \in E(\widetilde{H})$ with probability $\gl_{ij,uv}$, all
choices independent. The proofs of (\ref{inq-sec.mom.main}) and
(\ref{inq-sec.mom.main.2}) involve a second moment calculation. For
each $f \in \cH_g(G,H)$, set ${\bf X}_f = {\bf 1}_{\{f \in
\cH_g(G,\widetilde{H})\}}$ and ${\bf X}=\sum_{f \in \cH_g(G,H)} {\bf
X}_f$. Note that ${\bf X}=|\cH_g(G,\widetilde{H})|$. For each $f \in
\cH_g(G,H)$ we have
\begin{eqnarray}
{\mathbb E}({\bf X}_f) & = & {\mathbb P}(f \in
\cH_g(G,\widetilde{H}))
 \nonumber \\
 & = & {\mathbb P}(\{f(u)f(v)\in E(\widetilde{H})~\forall uv \in E(G)\})
\nonumber \\
& = & \prod_{uv\in E(G)}\gl_{g(u)g(v),uv}, \label{using.ind}
\end{eqnarray}
with (\ref{using.ind}) following from the fact that $\{f(u)f(v):uv
\in E(G)\}$ is a collection of disjoint edges and so $\{\{f(u)f(v)
\in E(\widetilde{H})\}:uv \in E(G)\}$ is a collection of independent
events. By linearity of expectation we therefore have
\begin{equation} \label{eq-mu}
{\mathbb E}({\bf X}) = |\cH_g(G,H)|\prod_{uv\in
E(G)}\gl_{g(u)g(v),uv} =C^Nw^W(g) := \mu.
\end{equation}

We now consider the second moment. For $f, f^\prime \in \cH_g(G,H)$
write $f \sim f^\prime$ if there is $uv \in E(G)$ with
$f(u)=f^\prime(u)$ and $f(v)=f^\prime(v)$. Note that ${\bf X}_f$ and
${\bf X}_{f^\prime}$ are not independent iff $f \sim f^\prime$. By
standard methods (see {\em e.g.} \cite{AlonSpencer}) we have
\begin{eqnarray}
\Var({\bf X}) & \leq & \mu + \sum_{(f,f^\prime) \in \cH_g(G,H)^2~:~f
\sim f^\prime} {\mathbb P}(\{f \in
\cH_g(G,\widetilde{H})\}\wedge\{f^\prime \in
\cH_g(G,\widetilde{H})\})
\nonumber \\
& \leq & \mu + |\{(f, f^\prime) \in \cH_g(G,H)^2:f \sim
f^\prime\}|. \nonumber
\end{eqnarray}
To estimate $|\{(f, f^\prime) \in \cH_g(G,H)^2:f \sim f^\prime\}|$
note that there are $|\cH_g(G,H)|$ choices for $f$, at most $N^2$
choices for a $uv \in E(G)$ on which $f$ and $f^\prime$ agree, and
finally at most
$$
\frac{|\cH_g(G,H)|}{C\gl_{g(u),u}C\gl_{g(v),v}} \leq
\frac{|\cH_g(G,H)|}{C^2 \gl_{vmin}^2}
$$
choices for the rest of $f^\prime$. We therefore have
\begin{eqnarray}
\frac{\Var({\bf X})}{\mu^2} & \leq & \frac{1}{\mu} +
\frac{|\cH_g(G,H)|^2N^2}{\mu^2C^2\gl_{vmin}^2} \nonumber \\
& = &
\frac{1}{C^2}\left(\frac{1}{C^{N-2}w^W(g)}+\frac{N^2}{\gl_{vmin}^2\prod_{uv
\in E(G)}\gl_{g(u)g(v),uv}}\right) \nonumber \\
& \leq & \frac{1}{C^2}\left(\frac{1}{w_{min}}+
\frac{\gl_{vmax}^N N^2}{\gl_{vmin}^2 w_{min}}\right) \label{exp.1} \\
& \leq & \frac{\alpha(N,a,b,W)}{C^2} \nonumber
\end{eqnarray}
for some function $\alpha$ (independent of $G$ and $g$). In
(\ref{exp.1}) we use (\ref{eq-mu}) and the fact that $N\geq 2$
(which holds since $G$ is non-empty). By Tchebychev's inequality, we
therefore have
$$
{\mathbb P}\left(|\mu - {\bf X}| > \mu\sqrt{\frac{\alpha}{C}}\right)
\leq \frac{1}{C}.
$$
It follows that the probability that $\widetilde{H}$ fails to
satisfy
\begin{equation} \label{last}
\left|\prod_{uv\in E(G)}\gl_{g(u)g(v),uv} -
\frac{|\cH_g(G,\widetilde{H})|}{|\cH_g(G,H)|}\right| \leq
\sqrt{\frac{\alpha}{C}}\prod_{uv\in E(G)}\gl_{g(u)g(v),uv}
\end{equation}
for a particular $g$ is at most $1/C$, and so the probability that
it fails to satisfy (\ref{last}) for any $g$ is at most $m^N/C$.

A similar argument gives that for a particular $v \in \cO_G$ and
$g:V(K_{a,b})\rightarrow \{1,\ldots,m\}$ the probability that
$\widetilde{H}$ fails to satisfy
\begin{equation} \label{last.v}
\left|\prod_{w_jz_k\in E(K_a,b)}\gl_{g(w_j)g(z_k),n_k(v)v} -
\frac{|\cH^v_g(K_{a,b},\widetilde{H})|}{|\cH^v_g(K_{a,b},H)|}\right|
\leq \sqrt{\frac{\alpha}{C}}\prod_{w_jz_k\in
E(K_a,b)}\gl_{g(w_j)g(z_k),n_k(v)v}
\end{equation}
is at most $1/C$, and so the probability that it fails to satisfy
(\ref{last.v}) for any $g$ is at most $m^{a+b}/C$.

As long as $C > m^N+aNm^{a+b}/(a+b)$ there is therefore an
$\widetilde{H}$ for which (\ref{inq-sec.mom.main}) is satisfied
for each $g:V(G)\rightarrow \{1,\ldots, m\}$ and
(\ref{inq-sec.mom.main.2}) is satisfied for each $v \in \cO_G$ and
$g:V(K_{a,b})\rightarrow \{1,\ldots, m\}$ with
$$
\delta(C)= \frac{\sqrt{\alpha}}{\sqrt{C}-\sqrt{\alpha}}.
$$
Since $\delta(C) \rightarrow 0$ as $C \rightarrow \infty$, we are
done.

\medskip

\noindent {\bf Acknowledgment} We are grateful to Alex Scott
\cite{Scott} for suggesting the construction of $\widetilde{H}$.


\begin{thebibliography}{99}

\bibitem{AlonSpencer} N. Alon and J. Spencer, {\em The
Probabilistic Method}, Wiley, New York, 2000.

\bibitem{CFGS} F.R.K. Chung, P. Frankl, R. Graham
and J.B. Shearer, Some intersection theorems for ordered sets and
graphs, {\em J. Combin. Theory Ser. A.} {\bf 48} (1986), 23--37.

\bibitem{Diestel}
R. Diestel, {\em Graph Theory}, Springer, New York, 1997.

\bibitem{GalvinTetali} D. Galvin and P. Tetali, On weighted graph homomorphisms, DIMACS Series in Discrete
Mathematics and Theoretical Computer Science {\bf 63} (2004) {\em
Graphs, Morphisms and Statistical Physics}, 97--104.

\bibitem{Kahn}
J. Kahn, An entropy approach to the hard-core model on bipartite
graphs, {\em Combin. Prob. Comp.} {\bf 10} (2001), 219--237.

\bibitem{McE} R.J. McEliece,
{\em The Theory of Information and Coding}, Addison-Wesley, London,
1977.

\bibitem{Scott} A. Scott, personal communication.

\end{thebibliography}
\end{document}